\documentclass[12pt]{amsproc}%{amsart}
\usepackage{amsmath,amssymb,amscd}
\newcommand{\cA}{{\mathcal A}}
\newcommand{\cB}{{\mathcal B}}
\newcommand{\cC}{{\mathcal C}}

\newcommand{\cF}{{\mathcal F}}
\newcommand{\cG}{{\mathcal G}}
\newcommand{\cH}{{\mathcal H}}

\newcommand{\cK}{{\mathcal K}}

\newcommand{\bbN}{{\mathbb N}}

\begin{document}

\title[Hybrid Normed Ideal Perturbations of $n$-tuples]{Hybrid Normed Ideal Perturbations of \\ $n$-tuples of Operators II: \\ Weak Wave Operators}
\author{Dan-Virgil Voiculescu}
\address{D.V. Voiculescu \\ Department of Mathematics \\ University of California at Berkeley \\ Berkeley, CA\ \ 94720-3840}
\thanks{Research supported in part by NSF Grant DMS-1665534.}
\dedicatory{Dedicated to the memory of Ronald G. Douglas.}
\date{}
\begin{abstract}
We prove a general weak existence theorem for wave operators for hybrid normed ideal perturbations. We then use this result to prove the invariance of Lebesgue absolutely continuous parts of $n$-tuples of commuting hermitian operators under hybrid normed ideal perturbations from a class studied in the first paper of this series.
\end{abstract}

\maketitle

\bigskip
\section{Introduction}
\label{sec1}

In \cite{4} we adapted basics of normed ideal perturbations to the hybrid setting and then turned to hybrid perturbations of $n$-tuples of commuting hermitian operators and found an unexpected result which also required a substantial amount of technical work. We showed that up to a factor of proportionality the modulus of quasicentral approximation with respect to the normed ideal $\cC^-_n$ is the same as the hybrid one with respect to $(\cC^-_{p_1},\dots,\cC^-_{p_n})$ when $p_1^{-1}+\dots+p_n^{-1} = 1$. In particular this implies that under such a hybrid perturbation the existence or the absence of $n$-dimensional Lebesgue absolutely continuous spectrum is preserved. To get in full generality that the Lebesgue absolutely continuous parts are actually preserved up to unitary equivalence requires some existence results for weak wave operators, which is the aim of the present paper (the particular case of integrable multiplicity functions could have been deduced directly from the formula for the modulus of quasicentral approximation).  Thus our aim here will be to extend one of the main results of \cite{3} showing that certain weak limits for the quantities which are considered in order to get wave operators, give rise to intertwiners with vanishing kernels and kernels of adjoints. In essence the extension is not far from the earlier result in \cite{3}, however the argument in \cite{3} is already rather intricate and having also to make a few technical improvements we felt the reader may not be too happy to get to fill in all these details as one of the so-called ``exercises left to the reader''. So we opted for a more detailed presentation of the proofs.

This paper has two more sections besides the introduction and references. Section~\ref{sec2} gives the hybrid existence result for weak wave operators. Section~\ref{sec3} is devoted to consequences, especially the invariance of the $n$-dimensional Lebesgue absolutely continuous parts.

This paper being the second one on hybrid normed ideal perturbations we use consistently the notation and definitions introduced in the first paper of the series.

\bigskip
\section{Existence of generalized wave operators}
\label{sec2}

The hybrid setting which we will use in this section involves a separable $C^*$-algebra $\cA$, $1 \in \cB \subset \cA$ a dense $*$-subalgebra with a countable basis as a vector space and $1 \in \cB_k \subset \cB$, $1 \le k \le n$, $*$-subalgebras of $\cB$, so that $\cB$ is generated by $\cup_{1 \le k \le n} \cB_k$ as an algebra. Let also $\varphi \in \cF([n])$.

If $\rho$ is a non-degenerate $*$-representation of $\cA$ on $\cH$, we define the $\varphi$-singular and $\varphi$-absolutely continuous projections $E_{\varphi}^0(\rho)$ and $E_{\varphi}(\rho)$ and the corresponding subspace $\cH_{\varphi}^0(\rho)$ and $\cH_{\varphi}(\rho)$ as follows. We consider all $p$-tuples, $p \in \bbN$, $\tau$ of operators in ${\tilde \beta} = \amalg_{1 \le j \le n} \cB_j$ and denote this by $\tau \subset {\tilde \beta}$ (all $p \in \bbN$ are considered). Then we form
\[
E_{\varphi}^0(\rho) = \bigwedge_{\tau \subset {\tilde \beta}} E_{\varphi_{\tau}}^0(\rho(\tau))
\]
where $\varphi_{\tau}(h) = \varphi(j)$ if $\tau(h) \in \cB_j$, $1 \le h \le p$, $1 \le j \le n$ (concerning $E_{\varphi}^0(\rho(\tau))$, $E_{\varphi}(\rho(\tau))$ see the definitions in section~6 of \cite{4}). Remark that
\[
E_{\varphi_{\tau_1 \amalg \tau_2}}^0(\rho(\tau_1 \amalg \tau_2)) \subset E_{\varphi_{\tau_1}}^0(\rho(\tau_1)) \wedge E_{\varphi_{\tau_2}}^0(\rho(\tau_2))
\]
so that for $E_{\varphi_{\tau}}(\rho(\tau)) = I - E_{\varphi_{\tau}}^0(\rho(\tau))$ we have
\[
E_{\varphi_{\tau_1 \amalg \tau_2}}(\rho(\tau_1 \amalg \tau_2)) \supset E_{\varphi_{\tau_1}}(\rho(\tau_1)) \vee E_{\varphi_{\tau_2}}(\rho(\tau_2)).
\]
We also define $E_{\varphi}(\rho) = I - E_{\varphi}^0(\rho)$ so that
\[
E_{\varphi}(\rho) = \bigvee_{\tau \subset {\tilde \beta}} E_{\varphi}(\rho(\tau)).
\]
We shall also use the notation $\cH_{\varphi}^0(\rho) = E_{\varphi}^0(\rho)\cH$, $\cH_{\varphi}(\rho) = E_{\varphi}(\rho)\cH$.

It is easy to infer the following extension of Proposition~6.1 \cite{4}.

\bigskip
\noindent
{\bf 2.1. Proposition.} {\em If $A_m = A_m^* \in \cK$, $m \in \bbN$ are so that
\[
\sup_{m \in \bbN} \|A_m\| < \infty \mbox{ and}
\]
\[
\lim_{m \to \infty} |[\rho(b),A_m]|_{\varphi(j)} = 0
\]
when $b \in \cB_j$, $1 \le j \le m$, then we have
\[
s - \lim_{m \to \infty} A_mE_{\varphi}(\rho) = 0.
\]
Moreover $\cH_{\varphi}^0(\rho)$ and $\cH_{\varphi}(\rho)$ are $\rho(\cA)$-invariant and the restrictions $\rho \mid \cH_{\varphi}^0(\rho)$ and $\rho \mid \cH_{\varphi}(\rho)$ are disjoint representations of $\cA$.
}

\bigskip
\noindent
{\bf {\em Proof.}} For all $p$-tuples $\tau \subset {\tilde \beta}$ by Prop.~6.1 \cite{4} we have that
\[
s - \lim_{m \to \infty} A_mE_{\varphi}(\rho(\tau)) = 0.
\]
Since the union of the $\cH_{\varphi(\tau)}(\rho(\tau)) = E_{\varphi(\tau)}(\rho(\tau))\cH$ is dense in $\cH_{\varphi}(\rho)$ we get that
\[
s - \lim_{m \to \infty} A_mE_{\varphi}(\rho) = 0.
\]
Clearly $\cH_{\varphi_{\tau}}^0(\rho(\tau)) = E_{\varphi_{\tau}}^0(\rho(\tau))\cH$ is invariant under $\rho(\tau)$ and hence their intersection over all $\tau \subset {\tilde \beta}$, $\cH_{\varphi}^0(\rho) = E_{\varphi}^0(\rho)\cH$ is invariant under $\rho({\tilde \beta})$, that is under $\rho(\cA)$. Also, since there are no non-zero $\rho(\tau)$-intertwiners between $E_{\varphi_{\tau}}^0(\rho(\tau))\cH$ and $E_{\varphi_{\tau}}(\rho(\tau))\cH$, there are no non-zero $\rho(\tau)$-intertwiners between $E_{\varphi}^0(\rho)\cH$ and $E_{\varphi_{\tau}}(\rho(\tau))\cH$. It follows also that there are no non-zero $\rho({\tilde \beta})$-intertwiners between $E_{\varphi}^0(\rho)\cH$ and $E_{\varphi}(\rho)\cH$.\qed

\bigskip
One of the facts which will be used in the proof of the main result of this section, a theorem which improves and extends Theorem~$1.4$ of \cite{3}, is a fact also used in the proof of the earlier result. If $\cG_{\Phi}^{(0)} \ne \cC_1$ then for every $X \in \cG_{\Phi}^{(0)}$ we have
\[
\lim_{j \to \infty} |\underset{j\mbox{-times}}{\underbrace{X \oplus \dots \oplus X}}|_{\Phi} = 0.
\]
When $X$ is a rank one orthogonal projection this is due to Kuroda (see \cite{1}, ch.~X, \S 2, the proof of Theorem~2.3). For general $X \in \cG_{\Phi}^{(0)}$ this then follows immediately from the fact that rank one projections are total in $\cG_{\Phi}^{(0)}$.

\bigskip
\noindent
{\bf 2.1. Theorem.} {\em Let $\varphi$ be such that $\varphi(j) \ne \Phi_1$, $1 \le j \le n$ and let $\rho_1,\rho_2$ be unital $*$-representations of $\cA$ on $\cH$ such that $\rho_1(b) - \rho_2(b) \in \cG_{\varphi(j)}^{(0)}$ if $b \in \cB_j$, $1 \le j \le n$. Assume moreover that there is a sequence of unitary elements $u_m \in Z(\cA)$, $m \in \bbN$, where $Z(\cA)$ is the center of $\cA$, such that
\[
w-\lim_{m \to \infty} \rho_1(u_m) = w-\lim_{m \to \infty} \rho_2(u_m) = 0
\]
and that the weak limit
\[
W = w-\lim_{m \to \infty} \rho_2(u_m^*)\rho_1(u_m)E_{\varphi}(\rho_1)
\]
exists. Then $W$ intertwines $\rho_1$ and $\rho_2$ and {\em Ker} $W = E_{\varphi}^0(\rho_1)$, {\em Ker} $W^* = E_{\varphi}^0(\rho_2)$. Moreover we have
\[
W^* = w-\lim_{m \to \infty} \rho_1(u^*_m)\rho_2(u_m)E_{\varphi}(\rho_2)
\]
and the representations $\rho_1|\cH_{\varphi}(\rho_1),\rho_2|\cH_{\varphi}(\rho_2)$ of $\cA$, are unitarily equivalent.
}

\bigskip
\noindent
{\bf {\em Proof.}} Since the $\cB_j$'s generate $\cB$, we have $\rho_2(b) - \rho_1(b) \in \cK$ for all $b \in \cB$ and hence since $\cB$ is dense in $\cA$ it follows that $\rho_1(a)-\rho_2(a) \in \cK$ for all $a \in \cA$. If $a \in \cA$ we have
\[
W\rho_1(a) - \rho_2(a)W = w-\lim_{m \to \infty} \rho_2(u^*_m)(\rho_1(a)-\rho_2(a))\rho_1(u_m)E_{\varphi}(\rho_1)
\]
and since $\rho_1(a)-\rho_2(a) \in \cK$ and $w-\lim_{m \to \infty} \rho_1(u_m) = 0$ we infer that
\[
s-\lim_{m \to \infty} \rho_2(u^*_m)(\rho_1(a)-\rho_2(a))\rho_1(u_m)E_{\varphi}(\rho_1) = 0
\]
and hence that $W\rho_1(a) = \rho_2(a)W$. This implies that $\rho_2\mid\overline{W\cH}$ is unitarily equivalent to a subrepresentation of $\rho_1\mid\overline{W^*\cH}$ and since $W = WE_{\varphi}(\rho_1)$ this is a subrepresentation of $\rho_1\mid\cH_{\varphi}(\rho_1)$ . We then must have $\overline{W\cH} \perp \cH_{\varphi}^0(\rho_2)$, that is $\overline{W\cH} \subset \cH_{\varphi}(\rho_2)$ or equivalently $W = E_{\varphi}(\rho_2)W$. Let ${\tilde W} = \cH_{\varphi}(\rho_2)\mid W \mid \cH_{\varphi}(\rho_1)$, that is the operator from $\cH_{\varphi}(\rho_1)$ to $\cH_{\varphi}(\rho_2)$ one gets from $W$. The main fact to be proved will be that Ker ${\tilde W} = 0$ and Ker ${\tilde W}^* = 0$. Before taking up this task, we shall prove a certain symmetry between $W$ and $W^*$. More precisely the symmetry is between $\rho_1,\rho_2,u_m,W$ and $\rho_2,\rho_1,u_m,W^*$. For this, we must show that the weak limit
\[
V = w - \lim_{m \to \infty} \rho_1(u_m^*)\rho_2(u_m)E_{\varphi}(\rho_2)
\]
exists and that $V = W^*$. Without assuming the existence of this weak limit, we can pass to a subsequence so that the weak limit defining $V$ exists and it will suffice to show that $V = W^*$ in this case, since then the operator $V$ we get will not depend on the chosen subsequence. Repeating for $V$ the argument with which we began the proof of the theorem, we find that $V$ intertwines $\rho_2$ and $\rho_1$ and hence that $E_{\varphi}(\rho_1)VE_{\varphi}(\rho_2) = V$. But then $V = w-\lim_{m \to \infty} E_{\varphi}(\rho_1)\rho_1(u^*_m)\rho_2(u_m)E_{\varphi}(\rho_2)$ and since $W = w-\lim_{m \to \infty} E_{\varphi}(\rho_2)\rho_2(u_m^*)\rho_1(u_m)E_{\varphi}(\rho)$ it follows that $V = W^*$.

\bigskip
To prove that Ker ${\tilde W} = 0$ and Ker ${\tilde W}^* = 0$, we shall assume the contrary and show that this leads to a contradiction. Let $P,Q$ be the orthogonal projections onto Ker ${\tilde W}$ and Ker ${\tilde W}^*$ respectively and remark that
\[
P \in ((\rho_1\mid\cH_{\varphi}(\rho_1))(\cA))',\ Q \in ((\rho_2\mid\cH_{\varphi}(\rho_2)))(\cA))'.
\]
In view of the symmetry we can assume that $P \ne 0$. Note also that rank $P$ must be infinite since otherwise $P \le E_{\varphi}^0(\rho_1)$.

The assumption $P \ne 0$ means there is $\xi \in \cH_{\varphi}(\rho_1)$, $\|\xi\| = 1$ so that
\[
w-\lim_{m \to \infty} W_m\xi = 0
\]
where $W_m = \rho_2(u^*_m)\rho_1(u_m)$. Replacing the $u_m$'s by a subsequence we may assume $p \ne q \Rightarrow \|W_p\xi - W_q\xi\| > 1$. Let then $A_k = W^*_{m_{k+1}}W_{m_k} - I$ for a sequence $m_1 < m_2 < \dots$ which we shall define recurrently. Since $P$ has infinite rank let $(\xi_k)_{k \in \bbN}$ be an orthonormal basis of $P\cH$ so that $\xi_1 = \xi$. Let further $\beta_j$ be a basis of the vector space $\cB_j$ and let $(b_r)_{r \in \bbN}$ be an enumeration of $\beta_1 \amalg \dots \amalg \beta_n$. In particular there is a map $\gamma: \bbN \to [n]$ so that $b_r \in \beta_{\gamma(r)}$ and $\beta_j = \{b_r \mid r \in \gamma^{-1}(j)\}$ for $1 \le j \le n$. We take $m_1 = 1$. Suppose $m_1 < \dots < m_k$ have been chosen. Then we can find $m_{k+1} > m_k$ so that
\[
\|(\rho_1(b_l)-\rho_2(b_l))\rho_2(u_{m_{k+1}})\rho_2(u^*_{m_k})\rho_1(u_{m_k})\xi_i\| < 1/k
\]
\[
\|(\rho_1(b_l)-\rho_2(b_l))\rho_2(u_{m_k})\rho_2(u^*_{m_{k+1}})\rho_1(u_{m_{k+1}})\xi_i\| < 1/k
\]
for $1 \le i$, $l \le k+1$. This is indeed possible because $\rho_1(b_l)-\rho_2(b_l) \in \cK$ and
\[
w-\lim_{m \to \infty}\rho_2(u_m) = 0,\ w-\lim_{m \to \infty}\rho_2(u^*_m)\rho_1(u_m)P = 0
\]
which implies that
\[
\lim_{m \to \infty}\|(\rho_1(b_l)-\rho_2(b_l))\rho_2(u_m)\rho_2(u^*_{m_k})\rho_1(u_{m_k})\xi_i\| = 0
\]
\[
\lim_{m \to \infty}\|(\rho_1(b_l)-\rho_2(b_l))\rho_2(u_{m_k})\rho_2(u^*_m)\rho_1(u_m)\xi_i\| = 0
\]
for all $i,l \in \bbN$.

For the above choice of the sequence $m_1 < m_2 < \dots$ we shall prove that
\[
s-\lim_{k \to \infty} [P A_kP,\rho_1(b)] = 0 \mbox{ and}
\]
\[
s-\lim_{k \to \infty} [P A^*_kP,\rho_1(b)] = 0
\]
for all $b \in \cB$. Since $\beta_1 \cup \dots \cup \beta_n$ generates $\cB$ as an algebra it will suffice to prove this when $b \in \beta_1 \cup \dots \cup \beta_n = \{b_r \mid r \in \bbN\}$.

We have $[P A_kP,\rho_1(b_r)] = P[A_k,\rho_1(b_r)]P$ and $[P A_k^*P,\rho_1(b_r)] = P[A_k^*,\rho_1(b_r)]P$ so that it will suffice to show that for all $i,r \in \bbN$ we have
\[
\begin{aligned}
\lim_{k \to \infty} \|[A_k,\rho_1(b_r)]\xi_i\| &= 0 \\
\lim_{k \to \infty} \|[A_k^*,\rho_1(b_r)]\xi_i\| &= 0.
\end{aligned}
\]
We have
\[
\begin{array}{rll}
[A_k,\rho_1(b_r)] &= &[W^*_{m_{k+1}}W_{m_k},\rho_1(b_r)] \\
&= &[\rho_1(u^*_{m_{k+1}})\rho_2(u_{m_{k+1}})\rho_2(u^*_{m_k})\rho_1(u_{m_k}),\rho_1(b_r)] \\
&= &\rho_1(u^*_{m_{k+1}})\rho_2(u_{m_{k+1}}u^*_{m_k})(\rho_1(b_r)-\rho_2(b_r))\rho_1(u_{m_k}) \\
&- &\rho_1(u^*_{m_{k+1}})(\rho_1(b_r)-\rho_2(b_r))\rho_2(u_{m_{k+1}})\rho_2(u^*_{m_k})\rho_1(u_{m_k}).
\end{array}
\]
This gives
\[
\begin{array}{l}
\displaystyle{\limsup_{k \to \infty}}\|[A_k,\rho_1(b_r)]\xi_i\| \\
\quad\le \displaystyle{\limsup_{k \to \infty}}(\|(\rho_1(b_r)-\rho_2(b_r))\rho_1(u_{m_k})\xi_i\| \\
\quad+ \|(\rho_1(b_r)-\rho_2(b_r))\rho_2(u_{m_{k+1}})\rho_2(u^*_{m_k})\rho_1(u_{m_k})\xi_i\|) = 0
\end{array}
\]
because $\rho_1(b_r) - \rho_2(b_r) \in \cK$, $w-\lim_{k\to\infty}\rho_1(u_{m_k}) = 0$ and because of the choice of $m_k$'s we made. Similarly we have
\[
\begin{array}{rll}
[A_k^*,\rho_1(b_r)] &= &\rho_1(u^*_{m_k})\rho_2(u_{m_k}u^*_{m_{k+1}})(\rho_1(b_r)-\rho_2(b_r))\rho_1(u_{m_{k+1}}) \\
&- &\rho_1(u_{m_k}^*)(\rho_1(b_r)-\rho_2(b_r))\rho_2(u_{m_k})\rho_2(u^*_{m_{k+1}})\rho_1(u_{m_{k+1}})
\end{array}
\]
so that
\[
\begin{array}{rll}
\displaystyle{\limsup_{k \to \infty}} \|[A_k^*,\rho_1(b_r)]\xi_i\| 
&\le &\displaystyle{\limsup_{k \to \infty}}(\|(\rho_1(b_r)-\rho_2(b_r))\rho_1(u_{m_{k+1}})\xi_i\| \\
&+ &\|(\rho_1(b_r)-\rho_2(b_r))\rho_2(u_{m_k})\rho_2(u^*_{m_{k+1}})\rho_1(u_{m_{k+1}})\xi_1\|) \\
&= &0
\end{array}
\]
again because $\rho_1(b_r)-\rho_2(b_r) \in \cK$, $w-\lim_{k\to\infty} \rho_1(u_{m_{k+1}}) = 0$ and of the choice of the $m_k$'s.

Remark now that $A_k + I$ being unitary we have
\[
A_k^*A_k = A_kA^*_k = -A_k - A_k^*.
\]
This then gives
\[
s - \lim_{k \to \infty}[P A_k^*A_kP,\rho_1(b_r)] = s - \lim_{k \to \infty} [\rho_1(b_r),P A_kP + P A_k^* P] = 0
\]
for all $r \in \bbN$. Since $\cB$ is self-adjoint this also gives
\[
s - \lim_{k \to \infty} ([P A_k^*A_k P,\rho_1(b_r)])^* = 0
\]
for all $r \in \bbN$. Further, since $u_n \in \cA$, we have $\rho_1(u_n)-\rho_2(u_n) \in \cK$ so that $W_n \in I + \cK$ and $A_k \in \cK$, $k \in \bbN$. The computations of $[A_k,\rho_1(b_r)]$ and $[A_k^*,\rho_1(b_r)]$ we did earlier in this proof, show that there are unitary operators $V_k,V'_k,V''_k,V'''_k,{\widetilde V}_k,{\widetilde V}'_k,{\widetilde V}''_k,{\widetilde V}'''_k$ so that
\[
[\rho_1(b_r),A_k] = V_k(\rho_1(b_r)-\rho_2(b_r))V'_k + V''_k(\rho_1(b_r)-\rho_2(b_r))V'''_k
\]
and
\[
[\rho_1(b_r),A_k^*] = {\widetilde V}_k(\rho_1(b_r)-\rho_2(b_r)){\widetilde V}'_k + {\widetilde V}''_k(\rho_1(b_r)-\rho_2(b_r)){\widetilde V}'''_k.
\]
It follows that
\[
\begin{array}{rll}
[PA^*_kA_kP,\rho_1(b_r)] &= &P[\rho_1(b_r),A_k+A_k^*]P \\
&= &P(V_k(\rho_1(b_r)-\rho_2(b_r))V'_k + V''_k(\rho_1(b_r)=\rho_2(b_r))V'''_k \\
&+ &{\widetilde V}'_k(\rho_1(b_r)-\rho_2(b_r)){\widetilde V}'_k + {\widetilde V}''_k(\rho_1(b_r)-\rho_2(b_r)){\widetilde V}'''_k)P \\
&\in &\cG^{(0)}_{\rho(\gamma(r))}.
\end{array}
\]
On the other hand the $*$-strong convergence of $[PA^*_kA_kP,\rho_1(b_r)]$ to $0$, easily gives that there are $k_1 < k_2 < \dots$ so that
\[
\begin{aligned}
&\lim_{j \to \infty} \left|j^{-1}|[P A_{k_1}^*A_{k_1}P,\rho_1(b_r)] + \dots + [P A_{k_j}^*A_{k_j}P,\rho_1(b_r)]|_{\varphi(\gamma(r))}\right. \\
&- j^{-1} \left.|([P A_{k_1}^*A_{k_1}P,\rho_1(b_r)]) \oplus \dots \oplus ([P A^*_{k_j}A_{k_j}P,\rho_1(b_r)])|_{\varphi(\gamma(r))}\right| \\
&= 0
\end{aligned}
\]
for all $r \in \bbN$.

In view of the result of the computation of $[P A_k^*A_kP,\rho_1(b_r)]$ we have
\[
\begin{aligned}
&|([P A_{k_1}^*A_{k_1}P_1,\rho_1(b_r)]) \oplus \dots \oplus ([P A_{k_j}^*A_{k_j}P,\rho_1(b_r)])|_{\varphi(\gamma(r))} \\
&\le 4|\underset{j\mbox{-times}}{\underbrace{(\rho_1(b_r)-\rho_2(b_r)) \oplus \dots \oplus (\rho_1(b_r)-\rho_2(b_r))}}|_{\varphi(\gamma(r))}.
\end{aligned}
\]
Since $\cG_{\varphi(\gamma(r))}^{(0)} \ne \cC_1$ we have
\[
\lim_{j \in \infty} j^{-1}|\underset{j\mbox{-times}}{\underbrace{(\rho_1(b_r)-\rho_2(b_r)) \oplus \dots \oplus (\rho_1(b_r)-\rho_2(b_r))}}|_{\varphi(\gamma(r))} = 0.
\]
Hence, if $B_j = j^{-1}(P A_{k_1}^*A_{k_1}P + \dots + P A_{k_j}^*A_{k_j}P)$ we have
\[
\lim_{j \to \infty} |[B_j,\rho_1(b_r)]|_{\varphi(\gamma(r))} = 0
\]
for all $r \in \bbN$, which then implies
\[
\lim_{j \to \infty} |[B_j,\rho_1(b)]|_{\varphi(\ell)} = 0
\]
for all $b \in \cB_{\ell}$.

Since $\|B_j\| \le 4$, $B_j \in \cK$ and $0 \le B_j \le P \le E_{\varphi}(\rho_1)$ it follows from Proposition~$2.1$ that
\[
s - \lim_{j \to \infty} B_j = 0.
\]
Recall now that $\xi_1 = \xi \in P\cH_{\varphi}(\rho_1)$ had the property that $p \ne q \Rightarrow \|W_p\xi - W_q\xi\| > 1$ which implies $\|A_kP\xi\| = \|A_k\xi\| = \|W_{m_{k+1}}^*W_{m_k}\xi - \xi\| > 1$ or equivalently $\langle P A_k^*A_kP\xi,\xi\rangle > 1$ for all $k \in \bbN$. This in turn implies $\langle B_j\xi,\xi\rangle > 1$, for all $j \in \bbN$ which is a contradiction.\qed

\bigskip
In the statement of Theorem~$2.1$ if we leave out the assumption that the weak limit
\[
w - \lim_{m \to \infty} \rho_2(u^*_m)\rho_1(u_m)E_{\varphi}(\rho_1)
\]
exists, it is always possible to find a subsequence of the $u_m$'s for which this weak limit exists and draw the conclusion that $\rho_1 \mid E_{\varphi}(\rho_1)$ and $\rho_2 \mid E_{\varphi}(\rho_2)$ are unitarily equivalent. Thus we have the following corollary.

\bigskip
\noindent
{\bf 2.1. Corollary.} {\em Let $\varphi$ be such that $\varphi(j) \ne \Phi_1$, $1 \le j \le n$ and let $\rho_1,\rho_2$ be unital $*$-representations of $\cA$ and $\cH$ such that $\rho_1(b)-\rho_2(b) \in \cG^{(0)}_{\varphi(j)}$ if $b \in \cB_j$, $1 \le j \le n$. Assume moreover that there is a sequence of unitary elements $u_m \in Z(\cA)$, $m \in \bbN$, where $Z(\cA)$ is the center of $\cA$, such that
\[
w - \lim_{m \to \infty} \rho_1(u_m) = w - \lim_{m \to \infty} \rho_2(u_m) = 0.
\]
Then the representations $\rho_1 \mid \cH_{\varphi}(\rho_1)$ and $\rho_2 \mid \cH_{\varphi}(\rho_2)$ of $\cA$ are unitarily equivalent.
}

\section{Invariance of Lebesgue absolutely continuous parts under perturbations}
\label{sec3}

\bigskip
\noindent
{\bf 3.1. Theorem.} {\em Let $\varphi \in \cF([n])$, $\varphi(j) = \Phi^-_{p_j}$, $p_j > 1$, $1 \le j \le n$, $n > 1$ be so that $p_1^{-1} + \dots + p_n^{-1} = 1$. Let $\tau$ and $\tau'$ be two $n$-tuples of commuting hermitian operators on $\cH$ so that $\tau(j) - \tau'(j) \in \cC_{p_j}^-$, $1 \le j \le n$. Then the Lebesgue absolutely continuous parts $\tau_{ac}$ and $\tau'_{ac}$ of $\tau$ and $\tau'$, are unitarily equivalent.
}

\bigskip
\noindent
{\bf {\em Proof.}} Consider the decompositions $\tau = \tau_{ac} \oplus \tau_s$, $\tau' = \tau'_{ac} \oplus \tau'_s$ with respect to $n$-dimensional Lebesgue measure and let $L > 0$ be such that $[-L,L]^n \supset \sigma(\tau) \cup \sigma(\tau')$. Recall also that by section~10 of \cite{4} these decompositions coincide with those into $\varphi$-singular and $\varphi$-absolutely continuous subspaces, in particular we have $k_{\varphi}(\tau_s) = k_{\varphi}(\tau'_s) = 0$. Consider also $\delta$ and $n$-tuple of multiplication operators by the coordinate functions in $L^2([-L,L]^n,d\lambda)$, where $\lambda$ is Lebesgue measure. If $\cA = C([-L,L]^n)$ is the $C^*$-algebra of continuous functions and $\cB \subset \cA$, the subalgebra of polynomial functions, with generator $\beta = (b_1,\dots,b_n)$ the $n$ coordinate functions we may form the representations of $\cA$ arising from functional calculus. Using Theorem~$5.1$ \cite{4} the adaptation of our non-commutative Weyl--von~Neumann type theorem, we find that $\tau_s \oplus \delta$ is unitarily equivalent mod the hybrid $n$-tuple $\cG_{\varphi}^{(0)}$ with $\delta$ and similarly $\tau'_s \oplus \delta$ is also unitarily equivalent with $\delta \mod \cG_{\varphi}^{(0)}$. This implies the existence of a unitary operator $U$ so that $\tau_{ac} \oplus \delta - U(\tau'_{ac} \oplus \delta)U^*$ is in $\cG_{\varphi}^{(0)}$. Let $\rho_1$ and $\rho_2$ be the representations of $\cA$ defined by $f \to f(\tau_{ac} \oplus \delta)$ and $f \to f(U(\tau'_{ac} \oplus \delta)U^*)$. Denoting by $\cB_j$ the subalgebra of $\cB$ consisting of polynomials in the $j$-th coordinate function we will have $\rho_1(b)-\rho_2(b) \in \cG_{\varphi(j)}^{(0)} = \cC_{p_j}^-$ if $b \in \cB_j$. Let further $u_m \in \cA$ be the function $u_m(x_1,\dots,x_n) = \exp(imx_1)$. Since the spectral measures of $\tau_{ac} \oplus \delta$ and $U(\tau'_{ac} \oplus \delta)U^*$ are absolutely continuous with respect to Lebesgue measure it is easily seen that $w - \lim_{m \to \infty} \rho_1(u_m) = w - \lim_{m \to \infty} \rho_2(u_m) = 0$. Thus the assumptions of Corollary~$2.1$ are satisfied and we get that $\rho_1$ and $\rho_2$ are unitarily equivalent (the singular parts being zero). This is in turn the same as the unitary equivalence of $\tau_{ac} \oplus \delta$ and $U(\tau'_{ac} \oplus \delta)U^*$ or $\tau'_{ac} \oplus \delta$. If $m_{ac}$ and $m'_{ac}$ are the multiplicity functions of $\tau_{ac}$ and $\tau'_{ac}$ we have proved that $m_{ac} + \chi_{[-L,L]^n}$ and $m'_{ac} + \chi_{[-L,L]^n}$ are equal almost everywhere with respect to Lebesgue measure. Clearly this implies $m_{ac} = m'_{ac}$ a.e.\ which is the unitary equivalence of $\tau_{ac}$ and $\tau'_{ac}$.\qed

\end{document}